# Loop-free Markov chains as determinantal point processes

## Alexei Borodin

*Department of Mathematics, California Institute of Technology, Pasadena, CA 91125, USA. E-mail: borodin@caltech.edu*

**Abstract.** We show that any loop-free Markov chain on a discrete space can be viewed as a determinantal point process. As an application, we prove central limit theorems for the number of particles in a window for renewal processes and Markov renewal processes with Bernoulli noise.

**Résumé.** Nous montrons que toute chaîne de Markov sans cycles sur un espace discret peut être vue comme un processus ponctuel determinantal. Comme application, nous démontrons des théorèmes limites centrales pour le nombre de particules dans une fenêtre pour des processus de renouvellement et des processus de renouvellement markoviens avec un bruit de Bernoulli.

## Introduction

Let $\mathfrak{X}$ be a discrete space. A (simple) random point process $\mathcal{P}$ on $\mathfrak{X}$ is a probability measure on the set $2^{\mathfrak{X}}$ of all subsets of $\mathfrak{X}$. $\mathcal{P}$ is called *determinantal* if there exists a $|\mathfrak{X}| \times |\mathfrak{X}|$ matrix $K$ with rows and columns marked by elements of $\mathfrak{X}$, such that for any finite $Y = (y_1, \ldots, y_n) \subset \mathfrak{X}$ one has

$$\mathcal{P}\{X \in 2^{\mathfrak{X}} | Y \subset X\} = \det[K_{y_i y_j}]_{i,j=1}^n.$$

The matrix $K$ is called a *correlation kernel* for $\mathcal{P}$.

A similar definition can be given for $\mathfrak{X}$ being any reasonable space; then the measure lives on locally finite subsets of $\mathfrak{X}$.

Determinantal point processes (with $\mathfrak{X} = \mathbb{R}$) have been used in random matrix theory since the early 60s. As a separate class, determinantal processes were first singled out in the mid-70s in [9] where the term *fermion* point processes was used. The term "determinantal" was introduced at the end of the 90s in [2], and it is now widely accepted. We refer the reader to surveys [1, 8, 12, 15] for further references and details.

Determinantal point processes have been enormously effective for studying scaling limits of interacting particle systems arising in a variety of domains of mathematics and mathematical physics including random matrices, representation theory, combinatorics, random growth models, etc. However, these processes are still considered as "exotic" – the common belief is that one needs a very special probabilistic model to observe a determinantal process.

The main goal of this note is to show that determinantal point processes are much more common. More exactly, we show that for *any* loop-free Markov chain the induced measure on trajectories is a determinantal point process. (Note that the absence of loops is essential – otherwise trajectories cannot be viewed as subsets of the phase space.) We work in a discrete state space in order to avoid technical difficulties like trajectories that may have almost closed loops, but our construction easily extends to suitable classes of loop-free Markov chains on continuous spaces as well.

Surprisingly, very little is known about arbitrary determinantal point processes. However, in the case when the correlation kernel is self-adjoint, several general results are available, and even more is known when the kernel is both self-adjoint and translation invariant (with $\mathfrak{X} = \mathbb{Z}^d$ or $\mathbb{R}^d$). See [1, 8, 12] and [15] for details.

In our situation, the kernel is usually *not* self-adjoint,[1] and self-adjoint examples should be viewed as "exotic." One such example goes back to [9], see also [12], Section 2.4: It is a 2-parameter family of renewal processes – processes on $\mathbb{Z}$ or $\mathbb{R}$ with positive i.i.d. increments. Our result implies that if we do not insist on self-adjointness then any process with positive i.i.d. increments is determinantal.

As an application, we prove a central limit theorem for the number of points in a growing window for renewal processes with Bernoulli noise and for Markov renewal processes (also known as semi-Markov processes) with Bernoulli noise. The proof is a version of the argument from [5, 13] adapted to non-self-adjoint kernels.

The key property that allows us to prove the central limit theorem is the boundedness of the operator defined by the correlation kernel in $\ell^2(\mathfrak{X})$. This boundedness is a corollary of what is known as a "renewal theorem" with a controlled rate of convergence.

## 1. Markov chains and determinantal point processes

Let $\mathfrak{X}$ be a discrete space, and let $P = [P_{xy}]_{x,y \in \mathfrak{X}}$ be the matrix of transition probabilities for a discrete time Markov chain on $\mathfrak{X}$. That is, $P_{xy} \geq 0$ for all $x, y \in \mathfrak{X}$ and

$$\sum_{y \in \mathfrak{X}} P_{xy} = 1 \quad \text{for any } x \in \mathfrak{X}.$$

Let us assume that our Markov chain is *loop-free*, i.e., the trajectories of the Markov chain do not pass through the same point twice almost surely. In other words, we assume that

$$(P^k)_{xx} = 0 \quad \text{for any } k > 0 \text{ and } x \in \mathfrak{X}.$$

This condition guarantees the finiteness of the matrix elements of the matrix

$$Q = P + P^2 + P^3 + \cdots.$$

Indeed, $(P^k)_{xy}$ is the probability that the trajectory started at $x$ is at $y$ after $k$th step. Hence, $Q_{xy}$ is the probability that the trajectory started at $x$ passes through $y \neq x$, and since there are no loops, we have $Q_{xy} \leq 1$. Clearly, $Q_{xx} \equiv 0$.

**Theorem 1.1.** *For any probability measure $\pi = [\pi_x]_{x \in \mathfrak{X}}$ on $\mathfrak{X}$, consider the Markov chain with initial distribution $\pi$ and transition matrix $P$ as a probability measure on trajectories viewed as subsets of $\mathfrak{X}$. Then this measure on $2^{\mathfrak{X}}$ is a determinantal point process on $\mathfrak{X}$ with correlation kernel*

$$K_{xy} = \pi_x + (\pi Q)_x - Q_{yx}.$$

---

[1] In fact, it can be written as a sum of a nilpotent matrix and a matrix of rank 1.

Before proceeding to the proof, let us point out that for any Markov chain $\{X(n)\}$ on a discrete state space $\widehat{\mathfrak{X}}$, its graph $\{(n, X(n)), n = 1, 2, \ldots\}$ is a loop-free Markov chain on $\mathfrak{X} = \mathbb{Z}_{>0} \times \widehat{\mathfrak{X}}$. Hence, the graph defines a determinantal point process.

Also, there exists a class of Markov chains $\{X(n)\}$ such that for each $n$, $X(n)$ is a random point configuration in some space $\widehat{\mathfrak{X}}_0$, and the graph $\{(n, X(n))\}$ defines a determinantal point process on $\mathbb{Z}_{>0} \times \widehat{\mathfrak{X}}_0$, see e.g. [3, 7, 16]. (Here $\widehat{\mathfrak{X}}$ is a suitable space of point configurations in $\widehat{\mathfrak{X}}_0$.) In such a case, the graph carries two types of determinantal structures, the second one on $\mathbb{Z}_{>0} \times \widehat{\mathfrak{X}}$ is afforded by Theorem 1.1.

The author is grateful to one of the referees for these remarks.

The proof of Theorem 1.1 is based on the following simple lemma, cf. [12], proof of Theorem 6.

**Lemma 1.2.** *For any* $a_1, \ldots, a_n, b_1, \ldots, b_n, c_1, \ldots, c_{n-1} \in \mathbb{C}$

$$\det \begin{bmatrix} a_1 b_1 & a_1 b_2 & a_1 b_3 & a_1 b_4 & \ldots & a_1 b_{n-1} & a_1 b_n \\ c_1 & a_2 b_2 & a_2 b_3 & a_2 b_4 & \ldots & a_2 b_{n-1} & a_2 b_n \\ * & c_2 & a_3 b_3 & a_3 b_4 & \ldots & a_3 b_{n-1} & a_3 b_n \\ * & * & c_3 & a_4 b_4 & \ldots & a_4 b_{n-1} & a_4 b_n \\ \ldots & \ldots & \ldots & \ldots & \ldots & \ldots & \ldots \\ * & * & * & * & \ldots & c_{n-1} & a_n b_n \end{bmatrix} = a_1(a_2 b_1 - c_1)(a_3 b_2 - c_2) \cdots (a_n b_{n-1} - c_{n-1})b_n.$$

**Proof.** The determinant in question is an $n$th degree polynomial in $a_i$'s with the highest coefficient $a_1 \cdots a_n b_1 \cdots b_n$, which vanishes at $a_1 = 0$, $a_2 = c_1/b_1$, $\ldots$, $a_n = c_{n-1}/b_{n-1}$. This implies the statement for $b_1 \cdots b_n \neq 0$, and thus, by continuity for all values of the parameters. □

**Proof of Theorem 1.1.** Let us evaluate the correlation function $\rho_n(x_1, \ldots, x_n)$ using the Markov chain. First, let us reorder the points $x_1, \ldots, x_n$ in such a way that $Q_{x_i x_j} = 0$ for $i > j$. This is always possible because if $Q_{xy} > 0$ then $Q_{yx} = 0$, and if $Q_{xy} > 0$ and $Q_{yz} > 0$ then $Q_{xz} > 0$, so we are simply listing the elements of a finite partially ordered set in a nondecreasing order. (For example, we can first list all minimal elements, then all elements that are minimal in the remaining set, etc.)

Once the points are ordered as described, using the Markov property we immediately compute

$$\rho_n(x_1, \ldots, x_n) = (\pi_{x_1} + (\pi Q)_{x_1}) Q_{x_1 x_2} \cdots Q_{x_{n-1} x_n}.$$

Lemma 1.2 (used with $b_i \equiv 1$) shows that this is exactly $\det[K_{x_i x_j}]_{i,j=1}^n$. □

Denote by $\mathfrak{Y}$ a finite subset of $\mathfrak{X}$ such that if a trajectory of our Markov chain leaves $\mathfrak{Y}$ then it does not return to $\mathfrak{Y}$ almost surely. In other words, if $Q_{yz} > 0$ for some $y \in \mathfrak{Y}$ and $z \in \mathfrak{X} \setminus \mathfrak{Y}$, then $Q_{zy'} = 0$ for any $y' \in \mathfrak{Y}$.

In order to consider the behavior of our Markov chain restricted to $\mathfrak{Y}$, it is convenient to contract the complement $\mathfrak{X} \setminus \mathfrak{Y}$ into one "final" state $\mathcal{F}$. Then the transition matrix for the new Markov chain on $\mathfrak{Y} \cup \{\mathcal{F}\}$ coincides with $P$ on $\mathfrak{Y} \times \mathfrak{Y}$, and for any $y \in \mathfrak{Y}$

$$P_{y\mathcal{F}} = 1 - \sum_{y' \in \mathfrak{Y}} P_{yy'} = \sum_{z \in \mathfrak{X} \setminus \mathfrak{Y}} P_{yz}, \qquad P_{\mathcal{F}y} = 0, \qquad P_{\mathcal{F}\mathcal{F}} = 1.$$

Set

$$\widetilde{\pi}_y(\mathfrak{Y}) = \text{Prob}\{\text{Intersection of the random trajectory with } \mathfrak{Y} \text{ starts at } y\},$$

$$\pi_0(\mathfrak{Y}) = \text{Prob}\{\text{The random trajectory does not intersect } \mathfrak{Y}\}.$$

Thus, $\pi_0(\mathfrak{Y}) = 1 - \sum_{y \in \mathfrak{Y}} \widetilde{\pi}_y(\mathfrak{Y})$. In what follows, we will use the notation $\widetilde{\pi}_y$ and $\pi_0$ instead of $\widetilde{\pi}_y(\mathfrak{Y})$ and $\pi_0(\mathfrak{Y})$ for the sake of brevity.

Let $K_{\mathfrak{Y}}$ denote the restriction of the matrix $K$ of Theorem 1.1 to $\mathfrak{Y} \times \mathfrak{Y}$.

**Theorem 1.3.** *In the above notation, we have* $\det(1 - K_{\mathfrak{Y}}) = \pi_0$. *If* $\pi_0 \neq 0$ *then the matrix* $L_{\mathfrak{Y}} := K_{\mathfrak{Y}}(1 - K_{\mathfrak{Y}})^{-1}$ *has the form*

$$(L_{\mathfrak{Y}})_{xy} = \frac{\widetilde{\pi}_x P_{y\mathcal{F}}}{\pi_0} - P_{yx}, \quad x, y \in \mathfrak{Y}.$$

**Proof.** The first statement is a well-known corollary of Theorem 1.1 which can also be easily proved using the inclusion–exclusion principle.

Assume that $\pi_0 \neq 0$. For any subset $Y$ of $\mathfrak{Y}$ choose the ordering $(y_1, \ldots, y_m)$ of its points in the same way as we did in the proof of Theorem 1.1: $Q_{y_i y_j} = 0$ for $i > j$. Then Lemma 1.2 yields

$$\det(L_{\mathfrak{Y}}|_{Y \times Y}) = \pi_0^{-1} \cdot \widetilde{\pi}_{y_1} P_{y_1 y_2} \cdots P_{y_{m-1} y_m} P_{y_m \mathcal{F}}$$

$$= \pi_0^{-1} \cdot \mathrm{Prob}\{\text{Intersection of the random trajectory with } \mathfrak{Y} \text{ is } Y\}.$$

Hence,

$$\det(1 + L_{\mathfrak{Y}}) = \sum_{Y \subset \mathfrak{Y}} \det(L_{\mathfrak{Y}}|_{Y \times Y}) = \pi_0^{-1},$$

and the matrix $(1 + L_{\mathfrak{Y}})$ is invertible.

Using the identity

$$\frac{L_{\mathfrak{Y}}}{1 + L_{\mathfrak{Y}}} = \left(\frac{-P_{\mathfrak{Y}}}{1 - P_{\mathfrak{Y}}}\right)^t + \frac{1}{1 - P_{\mathfrak{Y}}^t}(L_{\mathfrak{Y}} + P_{\mathfrak{Y}}^t)\frac{1}{1 + L_{\mathfrak{Y}}},$$

where by $P_{\mathfrak{Y}}$ we mean the restriction of $P$ to $\mathfrak{Y} \times \mathfrak{Y}$; and noting the fact that $P_{\mathfrak{Y}}(1 - P_{\mathfrak{Y}})^{-1}$ gives the restriction $Q_{\mathfrak{Y}}$ of $Q$ to $\mathfrak{Y} \times \mathfrak{Y}$ (here we use the assumption that the trajectories that leave $\mathfrak{Y}$ do not come back), we reduce the statement of the theorem to the following two relations:

$$\sum_{y' \in \mathfrak{Y}} \left(\frac{1}{1 - P_{\mathfrak{Y}}}\right)_{y'y} \widetilde{\pi}_{y'} = \pi_y + (\pi Q)_y, \qquad \sum_{y' \in \mathfrak{Y}} \left(\frac{1}{1 + L_{\mathfrak{Y}}}\right)_{y'y} P_{y'\mathcal{F}} \equiv \pi_0.$$

The two sides of the first relation represent two ways of computing the probability that $y \in \mathfrak{Y}$ lies on the random trajectory of our Markov chain, and the second relation is easily verified once both sides are multiplied by $(1 + L_{\mathfrak{Y}}^t)$. $\square$

## 2. Random point processes with Bernoulli noise

Let $\mathfrak{X}$ be a discrete set and let $\mathcal{P}$ be a random point process on $\mathfrak{X}$ (that is, $\mathcal{P}$ is a probability measure on $2^{\mathfrak{X}}$). Given two sequences of numbers

$$0 \leq p_x \leq 1, \qquad 0 \leq q_x \leq 1, \quad x \in \mathfrak{X},$$

we define a new random point process $\mathcal{P}^{(p,q)}$ on $\mathfrak{X}$ as follows. The random point configuration $X \subset \mathfrak{X}$ is being transformed by the following rules:

- Each particle $x \in X$ is removed with probability $p_x$, and it remains in place with probability $1 - p_x$;
- At each empty location $y \in \mathfrak{X} \setminus X$ a new particle appears with probability $q_y$, and the location remains empty with probability $1 - q_y$;
- The above operations at different locations are independent.

We say that $\mathcal{P}^{(p,q)}$ is the process $\mathcal{P}$ with *Bernoulli noise*.

Note that replacing the random point configuration $X$ with its complement on the whole space $\mathfrak{X}$ or on its part (the so-called particle–hole involution) can be viewed as a special case of the Bernoulli noise.

We will use the notation $D = \operatorname{diag}(d_x)$ for the diagonal matrix with matrix elements

$$D_{xy} = \begin{cases} d_x, & x = y, \\ 0, & x \neq y, \end{cases} \quad x, y \in \mathfrak{X}.$$

**Theorem 2.1.** *For a determinantal point process $\mathcal{P}$ on $\mathfrak{X}$ with correlation kernel $K$, the process with Bernoulli noise $\mathcal{P}^{(p,q)}$ is also determinantal, and its correlation kernel is given by*

$$K^{(p,q)} = \operatorname{diag}(q_x) + \operatorname{diag}(1 - p_x - q_x) \cdot K.$$

***Comments.*** *1. The theorem implies that if $p_x + q_x \equiv 1$ then the noisy process $\mathcal{P}^{(p,q)}$ is independent of the initial process $\mathcal{P}$, and it coincides with the product of Bernoulli random variables located at the points $x \in \mathfrak{X}$ with probability of finding a particle at $x$ equal to $q_x$. This fact is easy to see independently. Indeed, $1 - p_x$ and $q_x$ are conditional probabilities of finding a particle at $x$ given that the initial process has or does not have a particle there, and $p_x + q_x = 1$ implies that these probabilities are equal.*

*2. The particle–hole involution on $\mathfrak{X}$ or its part corresponds to taking both $p_x$ and $q_x$ equal to 1 on the corresponding part of the space.*

*3. The class of random point processes obtained as measures on trajectories of loop-free Markov chains is not stable under the application of Bernoulli noise.*

*4. Theorem 2.1 can be viewed as a special case of [4], Theorem 4.2.*

**Proof of Theorem 2.1.** We will consider the case when $p_x$ and $q_x$ are not equal to zero for just one $x = x_0 \in \mathfrak{X}$. The general case is clearly obtained by composition of such transformations.

From the definition of Bernoulli noise, we obtain

$$\mathcal{P}^{(p,q)}(X) = (1 - p_{x_0})\mathcal{P}(X) + q_{x_0}\mathcal{P}(X \setminus \{x_0\}) \quad \text{if } x_0 \in X.$$

Summing over all $X \in \mathfrak{X}$ that contain a fixed finite set $Y$ (which contains $x_0$) we obtain

$$\rho^{(p,q)}(Y) = (1 - p_{x_0}) \cdot \rho(Y) + q_{x_0} \cdot (\rho(Y \setminus \{x_0\}) - \rho(Y)) = (1 - p_{x_0} - q_{x_0}) \cdot \rho(Y) + q_{x_0} \cdot \rho(Y \setminus \{x_0\})$$

with

$$\rho(Y) = \mathcal{P}\{X \in 2^{\mathfrak{X}} | Y \subset X\}, \qquad \rho^{(p,q)}(Y) = \mathcal{P}^{(p,q)}\{X \in 2^{\mathfrak{X}} | Y \subset X\}.$$

One readily sees that the above relation is exactly the relation between the symmetric minors of the matrices $K^{(p,q)}$ and $K$ if the corresponding sets of rows and columns contain those marked by $x_0$. On the other hand, the correlation functions $\rho^{(p,q)}$ and $\rho$ of $\mathcal{P}^{(p,q)}$ and $\mathcal{P}$ away from $x_0$ are clearly identical, and same is true about the symmetric minors of $K^{(p,q)}$ and $K$ not containing the row and the column marked by $x_0$. □

## 3. CLT for number of points in a window of a determinantal process

Let $\mathfrak{X}$ be a discrete set and $\mathcal{P}$ be a determinantal process on $\mathfrak{X}$ with correlation kernel $K$. For any finite subset $\mathfrak{Y}$ of $\mathfrak{X}$, let us denote by $N_{\mathfrak{Y}}$ the number of points of the random point configuration that lie in $\mathfrak{Y}$. Using the definition of the correlation kernel, it is straightforward to verify that

$$\mathbb{E} N_{\mathfrak{Y}} = \operatorname{Tr} K_{\mathfrak{Y}}, \qquad \operatorname{Var} N_{\mathfrak{Y}} := \mathbb{E}(N_{\mathfrak{Y}} - \mathbb{E} N_{\mathfrak{Y}})^2 = \operatorname{Tr}(K_{\mathfrak{Y}} - K_{\mathfrak{Y}}^2),$$

where $\mathbb{E}$ denotes the expectation with respect to $\mathcal{P}$ and $K_{\mathfrak{Y}}$ is the restriction of $K$ to $\mathfrak{Y} \times \mathfrak{Y}$.

The following statement can be essentially extracted from [5, 13] and [14]. An important difference though is that here we consider correlation kernels which are not necessarily self-adjoint.

**Theorem 3.1.** *Assume that the correlation kernel $K$ defines a bounded operator in $\ell^2(\mathfrak{X})$. Then for any sequence $\{\mathfrak{Y}_m\}$ of finite subsets of $\mathfrak{X}$ such that*

$$\lim_{m\to\infty} |\mathfrak{Y}_m| = \infty, \qquad \operatorname{Var} N_{\mathfrak{Y}_m} \geq c_1 \cdot |\mathfrak{Y}_m|^{c_2}, \quad m \geq 1,$$

*for some strictly positive constants $c_1$ and $c_2$, the following central limit theorem for $N_{\mathfrak{Y}_m}$ holds:*

$$\lim_{m\to\infty} \operatorname{Prob}\left\{\frac{N_{\mathfrak{Y}_m} - \mathbb{E}N_{\mathfrak{Y}_m}}{\sqrt{\operatorname{Var} N_{\mathfrak{Y}_m}}} \leq x\right\} = \frac{1}{\sqrt{2\pi}} \int_{-\infty}^{x} e^{-t^2/2}\, dt$$

*for any $x \in \mathbb{R}$.*

**Proof.** As is shown in [13], Lemma 1, for $k \geq 2$ the $l$th cumulant of the random variable $N_{\mathfrak{Y}}$ is a linear combination of expressions of the form $\operatorname{Tr}(K_{\mathfrak{Y}} - K_{\mathfrak{Y}}^m)$ with $m$ ranging from 2 to $l$. Using the idea from [13], Lemma 2, we obtain

$$|\operatorname{Tr}(K_{\mathfrak{Y}} - K_{\mathfrak{Y}}^m)| \leq \sum_{k=1}^{m-1} |\operatorname{Tr}(K_{\mathfrak{Y}}^k - K_{\mathfrak{Y}}^{k+1})| \leq \sum_{k=1}^{m-1} \|K_{\mathfrak{Y}}^k - K_{\mathfrak{Y}}^{k+1}\|_1$$

$$\leq \|1 - K_{\mathfrak{Y}}\|_1 \sum_{k=1}^{m-1} \|K_{\mathfrak{Y}}\|^k \leq \|1 - K_{\mathfrak{Y}}\|_1 \sum_{k=1}^{m-1} \|K\|^k,$$

where we used the inequalities $|\operatorname{Tr} A| \leq \|A\|_1$ and $\|AB\|_1 \leq \|A\| \cdot \|B\|_1$. Since the trace norm of a matrix does not exceed the sum of the absolute values of its entries (and matrix elements of an operator are bounded by its norm), we see that

$$\|1 - K_{\mathfrak{Y}}\|_1 \leq (1 + \|K\|) \cdot |\mathfrak{Y}|^2.$$

Thus, the absolute value of any cumulant of $N_{\mathfrak{Y}_m}$ starting from the second one grows not faster than a fixed constant times $|\mathfrak{Y}_m|^2$. Hence, the absolute value of the $l$th cumulant of the normalized random variable $(N_{\mathfrak{Y}_m} - \mathbb{E}N_{\mathfrak{Y}_m})/\sqrt{\operatorname{Var} N_{\mathfrak{Y}_m}}$ is bounded by a constant times

$$\frac{|\mathfrak{Y}_m|^2}{(\operatorname{Var} N_{\mathfrak{Y}_m})^{l/2}} \leq c_1^{-1} \cdot |\mathfrak{Y}_m|^{2 - c_2 l/2}.$$

In particular, for large enough $l$ it converges to zero as $|\mathfrak{Y}_m| \to \infty$. [14], Lemma 3, completes the proof. $\square$

## 4. (Delayed) Renewal processes

In this section, we will apply our previous results to one of the simplest examples of loop-free Markov chains – renewal processes. For the sake of simplicity, we will consider only the discrete case. A large portion of our statements can be easily carried over to the continuous setting.

Let $\xi_0, \xi_1, \xi_2, \ldots$ be a sequence of mutually independent random variables with values in $\mathbb{Z}_{>0}$, such that $\xi_1, \xi_2, \ldots$ (but not $\xi_0$) have a common distribution. By definition, the sequence

$$S_0 = \xi_0, \qquad S_1 = \xi_0 + \xi_1, \qquad S_2 = \xi_0 + \xi_1 + \xi_2, \qquad \ldots$$

of partial sums forms a *delayed renewal process*, see, e.g., [6], Volume II, Section VI.6. The word "delayed" refers to a differently distributed $\xi_0$.

Clearly, $\{S_i\}_{i=0}^{\infty}$ is a trajectory of the loop-free Markov chain on $\mathbb{Z}_{>0}$ with the initial distribution and transition probabilities given by

$$\pi_i = \operatorname{Prob}\{\xi_0 = i\}, \qquad P_{ij} = \operatorname{Prob}\{\xi_1 = j - i\}, \quad i, j \in \mathbb{Z}_{>0}.$$

Thus, by Theorem 1.1, the sequence $\{S_i\}_{i=0}^\infty$ forms a determinantal point process on $\mathbb{Z}_{>0}$.

Let us now restrict our attention to the case when the distributions of random variables $\xi_i$ have exponentially decaying tails. It is very plausible that the results proved below can be extended to more general classes of distributions.

**Definition.** *We say that a random variable $\xi$ with values in $\mathbb{Z}_{>0}$ belongs to class $\mathcal{E}$ if there exists a constant $r \in (0,1)$ such that $\mathrm{Prob}\{\xi = n\} \leq r^n$ for large enough $n$.*

**Definition.** *An integral valued random variable $\xi$ is called* aperiodic *if*

$$\mathrm{Prob}\{\xi \in N\mathbb{Z}\} < 1 \quad \text{for any} \quad N \geq 2,$$

*that is, the distribution of $\xi$ is not fully supported by $N\mathbb{Z}$.*

The main result of this section is the following statement.

**Theorem 4.1.** *The delayed renewal process with aperiodic $\xi_1$ and $\xi_0, \xi_1 \in \mathcal{E}$, and with arbitrary Bernoulli noise defines a determinantal point process on $\mathbb{Z}_{>0}$ whose correlation kernel may be chosen so that it represents a bounded operator in $\ell^2(\mathbb{Z}_{>0})$.*

The proof is based on the following lemma.

**Lemma 4.2.** *Consider a random variable $\xi$ with values in $\mathbb{Z}_{>0}$ which belongs to $\mathcal{E}$ and is aperiodic. Set*

$$f(z) = \sum_{m=1}^\infty (\mathbb{E}z^\xi)^m = \sum_{n=1}^\infty f_n z^n.$$

*Then there exist positive constants $c_1$ and $c_2$ such that*

$$|f_n - (\mathbb{E}\xi)^{-1}| \leq c_1 e^{-c_2 n} \quad \text{for all } n \geq 1.$$

**Comments.** *The limit relation $\lim_{n\to\infty} f_n = (\mathbb{E}\xi)^{-1}$ holds without the assumption that $\xi \in \mathcal{E}$ and is commonly called the "Renewal theorem," see, e.g., [6], Volume I, Chapter XIII.*

**Proof of Lemma 4.2.** Denote $g(z) = \mathbb{E}z^\xi$. We have $|g(z)| < 1$ for $|z| < 1$, hence for $|z| < 1$

$$f(z) = \sum_{m=1}^\infty g^m(z) = \frac{g(z)}{1 - g(z)}.$$

Using the notation $g_1 = g'(1) = \mathbb{E}\xi$, we obtain

$$f(z) = \frac{1}{g_1(1-z)} + \frac{(1+g_1(1-z))g(z) - 1}{g_1(1-z)(1-g(z))}, \quad |z| < 1.$$

Since $\xi \in \mathcal{E}$, the function $g(z)$ is holomorphic in a disc of radius greater than 1, in particular, it is holomorphic in a neighborhood of 1. Hence, using the notation $g_2 = \frac{1}{2}g''(1) = \frac{1}{2}\mathbb{E}(\xi(\xi-1))$, we obtain

$$\frac{(1+g_1(1-z))g(z) - 1}{g_1(1-z)(1-g(z))} = -\frac{(1+g_1(1-z))(1+g_1(z-1)+g_2(z-1)^2 + \mathrm{O}((z-1)^3)) - 1}{g_1(1-z)^2(g_1 + \mathrm{O}(z-1))}$$

$$= \frac{g_1^2 - g_2 + \mathrm{O}(z-1)}{g_1^2(1 + \mathrm{O}(z-1))}$$

as $z \to 1$. Thus, $f(z) - 1/(g_1(1-z))$ is holomorphic in a neighborhood of $z = 1$. In addition to that, the fact that $\xi$ is aperiodic implies that $g(z)$ is not equal to 1 at any point $z \ne 1$ of the unit circle. Therefore, the function $f(z) - 1/(g_1(1-z))$ can be analytically continued to a disc of radius greater than 1. This implies the needed estimate. □

**Proof of Theorem 4.1.** By Theorem 2.1, it suffices to prove the statement without Bernoulli noise. Lemma 4.2 applied to $\xi = \xi_1$ implies that there exist positive constants $c_1$ and $c_2$ such that

$$|Q_{ij} - (\mathbb{E}\xi_1)^{-1}| \le c_1 e^{c_2(i-j)}, \quad i < j.$$

(Here we use the same notation as before: $Q = P + P^2 + P^3 + \cdots$.) We know that $\pi_i$'s decay exponentially as $i \to \infty$ because $\xi_0 \in \mathcal{E}$. Furthermore,

$$|(\pi Q)_j - (\mathbb{E}\xi)^{-1}| = \left| \sum_{i=1}^{j-1} \pi_i Q_{ij} - (\mathbb{E}\xi)^{-1} \right| \le (\mathbb{E}\xi)^{-1} \sum_{i=j}^{\infty} \pi_i + c_1 \sum_{i=1}^{j-1} \pi_i e^{c_2(i-j)},$$

which is easily seen to be exponentially decaying in $j \to \infty$.

Thus, the matrix

$$K_{ij} = \pi_i + (\pi Q)_i - Q_{ji}, \quad i,j = 1,2,\ldots,$$

afforded by Theorem 1.1 has the following properties:

$$|K_{ij}| \le \text{const}_1, \quad i \le j,$$
$$|K_{ij}| \le \text{const}_2 e^{-\text{const}_3(i-j)} + \text{const}_4 e^{-\text{const}_5 i}, \quad i > j.$$

for certain positive constants. This implies that if we consider the matrix $\widetilde{K}$ with matrix elements $\widetilde{K}_{ij} = e^{\alpha(i-j)} K_{ij}$ with positive $\alpha < \min\{\text{const}_3, \text{const}_5\}$ then the corresponding operator in $\ell^2(\mathbb{Z}_{\ge 0})$ will be bounded. On the other hand, the symmetric minors of $K$ and $\widetilde{K}$ are clearly the same. Thus, $\widetilde{K}$ is the needed correlation kernel. □

**Corollary 4.3.** *The central limit theorem for the number of particles in a window as described in Theorem 3.1 holds for the delayed renewal processes with Bernoulli noise which satisfy the assumptions of Theorem 4.1.*

**Comments.** *1. In the absence of Bernoulli noise, the central limit theorem for $N_{[0,T]}$ is an easy corollary of the classical central limit theorem for sums of i.i.d. random variables, see, e.g., [6], Volume II, Section XI.5. However, with the presence of Bernoulli noise the statement does not look as obvious.*

*2. In concrete examples (including the one mentioned in the previous comment), the asymptotic behavior of the expectation and variance of $N_\mathfrak{Y}$ can often be computed explicitly using the formulas*

$$\mathbb{E} N_\mathfrak{Y} = \text{Tr } K_\mathfrak{Y}, \qquad \text{Var } N_\mathfrak{Y} = \text{Tr}(K_\mathfrak{Y} - K_\mathfrak{Y}^2).$$

## 5. Markov renewal processes

Markov renewal process (or semi-Markov processes) are hybrids of Markov chains and renewal processes introduced by Lévy and Smith in 1954. An excellent introduction to Markov renewal processes can be found in [10] and [11].

We will use one of the simplest possible settings.

Consider a Markov chain on $\mathfrak{X} = S \times \mathbb{Z}_{>0}$ with transition probabilities

$$P((s_1, t_1) \to (s_2, t_2)) = P_{s_1 s_2}(t_2 - t_1),$$

where for any $s_1, s_2 \in S$ the function $P_{s_1 s_2}(t)$ is supported by $\mathbb{Z}_{>0}$,

$$\sum_{t=1}^{\infty} P_{s_1 s_2}(t) =: P_{s_1 s_2} \in [0,1], \quad s_1, s_2 \in S,$$

and $[P_{s_1 s_2}]_{s_1, s_2 \in S}$ is the matrix of transition probabilities for a Markov chain on $S$.

The projection of thus defined Markov chain on $\mathfrak{X}$ to the first coordinate is the Markov chain on $S$ sometimes called the "driving" Markov chain. The Markov chain on $\mathfrak{X}$ is a *Markov renewal process*. It can be viewed as the driving Markov chain with randomly transformed time scale – the passage time from $s_1$ to $s_2$ is a random variable depending on $s_1$ and $s_2$.

**Theorem 5.1.** *Assume that $S$ is finite and the driving Markov chain on $S$ is irreducible. Assume further that whenever $P_{ss'} \neq 0$, the passage time between $s$ and $s'$ (which is equal to $t$ with probability $P_{ss'}(t)/P_{ss'}$) is aperiodic and is in $\mathcal{E}$, and that the initial distribution for the Markov renewal process conditioned on the first coordinate being fixed (and arbitrary) is also in $\mathcal{E}$.*

*Then the point process on $\mathfrak{X}$ formed by trajectories of the Markov renewal process, with any Bernoulli noise, is a determinantal point process, and the correlation kernel may be chosen so that it represents a bounded operator in $\ell^2(\mathfrak{X})$.*

***Comments.*** 1. The assumptions of $S$ being finite, the driving chain being irreducible and passage times being aperiodic and in $\mathcal{E}$ can probably be relaxed. For example, instead of requiring that all passage times are aperiodic, one can impose the weaker condition of the first return time from any state to itself being aperiodic.

2. *Theorem 4.1 is a special case of Theorem 5.1 obtained by considering the set $S$ with a single element.*

**Proof of Theorem 5.1.** If we start the Markov renewal process at any state $s$ then the times of its arrivals to any state $s'$ form a delayed renewal process (see the previous section) with the initial delay $\xi_0$ being the first passage time from $s$ to $s'$ and the common distribution of $\xi_1, \xi_2, \ldots$ being the first return time from $s'$ to $s'$.

Both $\xi_0$ and $\xi_1$ are in $\mathcal{E}$; see Lemma 5.3. Furthermore, the aperiodicity of $\xi_1$ immediately follows from the assumption of all passage times being aperiodic. Since the initial distribution of our Markov renewal process conditioned on fixing the first coordinate is also in $\mathcal{E}$, when we start from this initial distribution the arrival times to any state also form a delayed renewal process with same renewal time $\xi_1$, but generally speaking, different initial delay time $\xi_0$, which however is still in $\mathcal{E}$.

Thus, we can apply the arguments used in the proof of Theorem 4.1 to produce an exponential bound on the difference of $Q((s, t_1) \to (s', t_2))$ and a known constant (equal to the inverse of the expectation of the first return time from $s'$ to $s'$).

The remainder of the proof is just the same as in Theorem 4.1: We conjugate the kernel afforded by Theorem 1.1 by a suitable diagonal matrix whose nonzero matrix elements form a geometric progression, and thus get the desired correlation kernel which defines a bounded operator in $\ell^2(\mathfrak{X})$. □

**Corollary 5.2.** *The central limit theorem for the number of particles in a window as described in Theorem 3.1 holds for the Markov renewal processes with Bernoulli noise which satisfy the assumptions of Theorem 5.1.*

It remains to prove the following lemma.

**Lemma 5.3.** *Under the assumptions of Theorem 5.1, for any $s, s' \in S$ the first passage time from $s$ to $s'$ is a random variable of class $\mathcal{E}$.*

**Proof.** For any $s_1, s_2 \in S$ such that $P_{s_1 s_2} \neq 0$ set

$$f_{s_1 s_2}(z) = \frac{1}{P_{s_1 s_2}} \sum_{t=1}^{\infty} P_{s_1 s_2}(t) z^t.$$

By the hypothesis, all such functions have analytic continuation to a disc of radius larger than some constant $R > 1$. Denote

$$M = \max_{|z|=R, \, s_1, s_2 \in S} |f_{s_1 s_2}(z)|.$$

Observe that $M \geq 1$ by the maximum principle because $f_{s_1 s_2}(1) \equiv 1$.

The probability of the Markov renewal process started at $s_1$ to walk the path $s_1 \to s_2 \to \cdots \to s_m$ and spend time $T$ on this path can be estimated as follows:

$$\sum_{t_1 + \cdots + t_{m-1} = T} P_{s_1 s_2}(t_1) P_{s_2 s_3}(t_2) \cdots P_{s_{m-1} s_m}(t_{m-1}) = P_{s_1 s_2} \cdots P_{s_{m-1} s_m} \cdot \frac{1}{2\pi \mathrm{i}} \int_{|z|=R} \frac{f_{s_1 s_2}(z) \cdots f_{s_{m-1} s_m}(z) \, \mathrm{d}z}{z^{T+1}}$$

$$\leq P_{s_1 s_2} \cdots P_{s_{m-1} s_m} \cdot M^{m-1}/R^T.$$

Note that $P_{s_1 s_2} \cdots P_{s_{m-1} s_m}$ is exactly the probability of the driving Markov chain started at $s_1$ to walk the path $s_1 \to \cdots \to s_m$.

Let $p$ be the minimum of all $P_{s_1 s_2} \neq 0$. Since the driving chain is irreducible, for any initial distribution its probability of hitting any given state in the first $|S|$ steps is at least $p^{|S|}$. Hence, the number of steps in the first passage of the driving chain from any $s_1 \in S$ to any $s_2 \in S$ is a random variable from $\mathcal{E}$. (Indeed, the probability of this number of steps being at least $n$ is $\leq (1 - p^{|S|})^{[n/|S|]}$.)

Take any $a > 0$ such that $M^a < R$. For a first passage of the Markov renewal process from $s_1$ to $s_2$ that takes time $T$, either the number of steps of the driving Markov chain is at least $[aT]$ or it is $< [aT]$. The probability of the former event decays exponentially in $T$ by the previous paragraph, while the probability of the latter event, by the estimate above, does not exceed $M^{[aT]}/R^T$, which also decays exponentially in $T$ by the choice of $a$. $\square$

## Acknowledgments

This research was partially supported by the NSF grant DMS-0402047 and the CRDF grant RIM1-2622-ST-04. The author would also like to thank the referees for a number of valuable suggestions.